\documentclass[11]{article}
\usepackage[russian]{babel}
\usepackage{amssymb}
\usepackage{inputenc}
\usepackage{graphicx}
\begin{document}
\begin{center}
\textbf{\Large PROPERTIES OF DISCRETE ANALOGUE OF THE DIFFERENTIAL
OPERATOR  $\frac{d^{2m}}{dx^{2m}}-\frac{d^{2m-2}}{dx^{2m-2}}$ }\\
\textbf{\large Kh.M.Shadimetov, A.R.Hayotov}
\end{center}

\textbf{Abstract.} In the paper properties of the discrete
analogue $D_m(h\beta)$ of the differential operator
$\frac{d^{2m}}{dx^{2m}}-\frac{d^{2m-2}}{dx^{2m-2}}$ are studied.
It is known, that zeros of differential operator
$\frac{d^{2m}}{dx^{2m}}-\frac{d^{2m-2}}{dx^{2m-2}}$ are functions
$e^x$, $e^{-x}$ and $P_{2m-3}(x)$. It is proved that discrete
analogue $D_m(h\beta)$ of this differential operator also have
similar properties.

\large

Present work is dedicated to investigation of properties of the
discrete analogue of one differential operator. We give necessary
definitions and formulas.

{\bf Definition 1.} {\it Function
$\varphi[\beta]=\varphi[\beta_1,\beta_2,...,\beta_n]$ is called by
function of discrete argument, if it given on some set of integer
values $\beta=(\beta_1,\beta_2,...,\beta_n)$.}

{\bf Definition 2.} {\it By inner product $[\varphi,\psi]$ of two
discrete functions $\varphi[\beta]$ and $\psi[\beta]$ is called
the number
$$
[\varphi,\psi]=\sum\limits_{\beta\in
B}\varphi[\beta]\bar\psi[\beta], \eqno (1)
$$
if the series in right hand side of (1) converge absolutely. Here
$B$ is a domain of definition of $\varphi[\beta]$ and
$\psi[\beta]$. }

{\bf Definition 3.}  {\it By convolution
$\varphi[\beta]*\psi[\beta]$ of two functions $\varphi[\beta]$ and
$\psi[\beta]$ is called the inner product}
$$
\chi[\beta]=\varphi[\beta]*\psi[\beta]=[\varphi[\gamma],\bar\psi[\beta-\gamma]].
$$

Following formula holds [1]
$$
\sum\limits_{\gamma=0}^{n-1}q^k\gamma^k=\frac{1}{1-q}\sum_{i=0}^k\left(\frac{q}{1-q}\right)
^i\Delta^i0^k-\frac{q^n}{1-q}\sum_{i=0}^k\left(\frac{q}{1-q}\right)
^i\Delta^i\gamma^k|_{\gamma=n}. \eqno (2) $$ When $|q|<1$ from (2)
we have
$$
\sum\limits_{\gamma=0}^{\infty}q^k\gamma^k=\frac{1}{1-q}\sum_{i=0}^k\left(\frac{q}{1-q}\right)
^i\Delta^i0^k. \eqno (3)
$$

Euler polynomial $E_k(\lambda)$ [2] is given by following
formula
$$
E_k(\lambda)=\frac{(1-\lambda)^{k+2}}{\lambda}D^k\frac{\lambda}{(1-\lambda)^2},
\mbox{ where } D=\lambda\frac{d}{d\lambda},\
D^k=\lambda\frac{d}{d\lambda} D^{k-1}. \eqno (4)
$$

For Euler polynomials following theorem is valid [3]

{\bf Theorem 1.}  {\it Polynomials
$$
P_k(x)=(x-1)^{k+1}\sum\limits_{i=1}^{k+1}\frac{\Delta^i0^{k+1}}{(x-1)^i}
$$
and
$$
P_k(\frac{1}{x})=(\frac{1}{x}-1)^{k+1}\sum\limits_{i=1}^{k+1}\left(\frac{x}{1-x}\right)\Delta^i0^{k+1}
$$
are Euler polynomials $E_k(x)=P_k(x)$ and
$E_k(\frac{1}{x})=P_k(\frac{1}{x})$ respectively. }

In construction of optimal quadrature and cubature formulas the
discrete analog $D_{m,n}[\beta]$ of the polyharmonic  operator
$\Delta^m=(\frac{\partial^2}{\partial x_1^2}+
\frac{\partial^2}{\partial x_2^2}+...+\frac{\partial^2}{\partial
x_n^2})^m$ plays main role in $L_2^{(m)}(R^n)$ space. First in
$L_2^{(m)}(R^n)$ space construction and investigation of
properties of convolution of discrete operator with function
$G_{m,n}[\beta]$, where $G_{m,n}(x)$ is a fundamental solution of
the polyharmonic operator, were studied by Sobolev [4]. The
discrete argument function $D_{m,n}[\beta]$ satisfies following
equality
$$
D_{m,n}[\beta]*G_{m,n}[\beta]=\delta[\beta],
$$
where $\delta[\beta]$ is equal to 1 as $\beta=0$, is equal to 0 as
$\beta\neq 0$. Sobolev gave the algorithm for finding of function
$D_{m,n}[\beta]$ and proved several properties of this function.

In one dimensional case, i.e. in $L_2^{(m)}(R)$ space discrete
analogue of the differential operator $\frac{d^{2m}}{dx^{2m}}$ was
constructed by Z.J.Jamalov [5,6]. But there were not found $m+1$
unknown coefficients. In the work [7] were found expressions of
that unknown coefficients and discrete analogue of the
differential operator $\frac{d^{2m}}{dx^{2m}}$ was completely
constructed.

In construction of the optimal quadrature formulas in
$W_2^{(m,m-1)}(0,1)$ space with the help of Sobolev's algorithm
[4] as above we need to construct a discrete analogue  of the
differential operator
$\frac{d^{2m}}{dx^{2m}}-\frac{d^{2m-2}}{dx^{2m-2}}$ instead of
operator $\frac{d^{2m}}{dx^{2m}}$. Here $W_2^{(m,m-1)}$ is Hilbert
space and norm of a function from this space is given by formula:
$$
\|\varphi(x)|W_2^{(m,m-1)}(0,1)\|=\left[\int\limits_0^1\left(
\varphi^{(m)}(x)+\varphi^{(m-1)}(x)\right)^2dx\right]^{1/2}.
$$

In work [8] was constructed the discrete analogue $D_m(h\beta)$ of
the differential operator
$\frac{d^{2m}}{dx^{2m}}-\frac{d^{2m-2}}{dx^{2m-2}}$, which
satisfies following equality
$$
D_m(h\beta)*\psi_m(h\beta)=\delta(h\beta), \eqno (5)
$$
where
$$
\psi_m(x)=\frac{\mathrm{sign}
x}{2}\left(\frac{e^x-e^{-x}}{2}-\sum\limits_{k=1}^{m-1}\frac{x^{2k-1}}
{(2k-1)!}\right). \eqno (6)
$$
is a solution of equation
$$
\left(\frac{d^{2m}}{dx^{2m}}-\frac{d^{2m-2}}{dx^{2m-2}}\right)\psi(x)=\delta(x),
\eqno (7)
$$
$$
\delta(h\beta)= \left\{ \begin{array}{ll} 1,& \beta=0,\\
0,&\beta\neq 0.
\end{array}
\right. \eqno (8) $$ Discrete function $D_m(h\beta)$ plays
important role in construction of optimal quadrature formulas in
$W_2^{(m,m-1)}(0,1)$ space. For the discrete operator
$D_m(h\beta)$ following theorem holds.

{\bf Theorem 2} [8]. {\it Discrete analogue of the differential
operator $\frac{d^{2m}}{dx^{2m}}-\frac{d^{2m-2}}{dx^{2m-2}}$
satisfying equality  (5) has following form:
$$
D_m(h\beta)=\frac{1}{p_{2m-2}^{(2m-2)}}\left\{
\begin{array}{ll}
\sum\limits_{k=1}^{m-1}A_k\lambda_k^{|\beta|-1}, & |\beta|\geq 2;\\
-2e^h+\sum\limits_{k=1}^{m-1}A_k, &|\beta|=1;\\
2C+\sum\limits_{k=1}^{m-1}\frac{A_k}{\lambda_k},& \beta=0,
\end{array}
\right. \eqno (9)
$$
where $$ C=1+(2m-2)e^h+e^{2h}+\frac{e^h\cdot
p_{2m-3}^{(2m-2)}}{p_{2m-2}^{(2m-2)}}, \eqno (10)
$$
$$
A_k=\frac{2(1-\lambda_k)^{2m-2}[\lambda_k(e^{2h}+1)-e^h(\lambda_k^2+1)]
p_{2m-2}^{(2m-2)}}{\lambda_kP_{2m-2}'(\lambda_k)}, \eqno (11)$$
$$
P_{2m-2}(\lambda)=\sum\limits_{s=0}^{2m-2}p_s^{(2m-2)}\lambda^s=(1-e^{2h})
(1-\lambda)^{2m-2}-2(\lambda(e^{2h}+1)-e^h(\lambda^2+1))\times
$$
$$ \times
\left[h(1-\lambda)^{2m-4}+\frac{h^3(1-\lambda)^{2m-6}}{3!}E_2(\lambda)+
...+\frac{h^{2m-3}E_{2m-4}(\lambda)} {(2m-3)!}\right], \eqno
(12)$$ $p_{2m-2}^{(2m-2)}$, $p_{2m-3}^{(2m-2)}$ are coefficients
of polynomial $P_{2m-2}(\lambda)$, $\lambda_k$ are roots of
polynomial $P_{2m-2}(\lambda)$, $|\lambda_k|<1$ , $E_k(\lambda)$
is Euler polynomial [2].}

Aim of the present work is investigation of properties of the
discrete analogue $D_m(h\beta)$ of the differential operator
$\frac{d^{2m}}{dx^{2m}}-\frac{d^{2m-2}}{dx^{2m-2}}$. It is known,
that zeros of the differential operator
$\frac{d^{2m}}{dx^{2m}}-\frac{d^{2m-2}}{dx^{2m-2}}$ are $e^x$,
$e^{-x}$ functions and a polynomial $P_{2m-3}(x)$ of  degree
$2m-3$. The discrete operator $D_m(h\beta)$ also has similar
properties, i.e. following is valid

\textbf{Theorem 3.} \emph{The discrete analogue $D_m(h\beta)$ of
the differential operator
$\frac{d^{2m}}{dx^{2m}}-\frac{d^{2m-2}}{dx^{2m-2}}$ when $m=1,2,3$
satisfies following equalities}

1) $D_m(h\beta)*e^{h\beta}=0;$

2) $D_m(h\beta)*e^{-h\beta}=0;$

3) \emph{$D_m(h\beta)*(h\beta)^n=0,$ $n\leq 2m-3$, i.e.
convolution of $D_m(h\beta)$ with a polynomial of degree $\leq
2m-3$ when $m=2,3$ is equal to zero};

4) \emph{$D_m(h\beta)*\psi_m(h\beta)=\delta(h\beta),$ where
$\psi_m(h\beta)$ and $\delta(h\beta)$ are defined by (6) and (8),
respectively.}
\newpage

\textbf{Proof of theorem 3.}

We prove the theorem by direct calculation of all convolutions.

\textbf{1).} For convenience leading coefficient of the polynomial
$P_{2m-2}(\lambda)$ we denote by:
$$
p=p_0^{(2m-2)}=p_{2m-2}^{(2m-2)}=(1-e^{2h})+2e^h\left(h+\frac{h^3}{3!}+...+
\frac{h^{2m-3}}{(2m-3)!}\right). \eqno(13) $$ Using the definition
of convolution of discrete functions and equality (9), we
calculate convolution $D_m(h\beta)*e^{h\beta}$:
$$
F_1=D_m(h\beta)*e^{h\beta}=\sum\limits_{\gamma=-\infty}^{\infty}D_m(h\gamma)e^{h\beta-h\gamma}=
e^{h\beta}\sum\limits_{\gamma=-\infty}^{\infty}D_m(h\gamma)e^{-h\gamma}=
$$
$$
=\frac{e^{h\beta}}{p}\left(\sum_{k=1}^{m-1}\frac{A_k}{\lambda_k}\sum_{\gamma=1}^{\infty}
\left[(\lambda_ke^h)^\gamma+(\frac{\lambda_k}{e^h})^\gamma\right]-2(e^{2h}+1)+2C+\sum_{k=1}^{m-1}
\frac{A_k}{\lambda_k}\right). \eqno(14)
$$
Since $|\lambda_k e^h|<1$ and $|\frac{\lambda_k}{e^h}|<1$, then
infinite series in equality (14) are convergent. Then
$$
F_1=
\frac{e^{h\beta}}{p}\left(\sum_{k=1}^{m-1}\frac{A_k}{\lambda_k}
\left[\frac{\lambda_ke^h}{1-\lambda_ke^h}+
\frac{\lambda_k}{e^h-\lambda_k}\right]-2(e^{2h}+1)+2C+\sum_{k=1}^{m-1}
\frac{A_k}{\lambda_k}\right).
$$
Hence, taking into account equalities (13), (11) and after some
simplifications, we get
$$
F_1=\frac{2e^{h\beta+h}}{p}\sum_{k=1}^{m-1}\left(\frac{(1-\lambda_k)^{2m-2}}{\lambda_k
\prod\limits_{i=1,i\neq k\atop i\neq
2m-1-k}^{2m-2}(\lambda_k-\lambda_i)}-
\frac{(1-\lambda_k)^2}{\lambda_k}\right). \eqno (15) $$ Let $m=1$,
then from (15) we obtain, that
$F_1=D_1(h\beta)*e^{h\beta}=0$.\\
Let $m=2$, then equality (15) gives
$$
F_1=D_2(h\beta)*e^{h\beta}=\frac{2e^{h\beta+h}}{p}\left(\frac{(1-\lambda_k)^2}{\lambda_k}-
\frac{(1-\lambda_k)^2}{\lambda_k}\right)=0
$$
Now, let $m=3$, then from (15) we get
$$
F_1=D_3(h\beta)*e^{h\beta}=
\frac{2e^{h\beta+h}}{p}\sum_{k=1}^{m-1}\left(\frac{(1-\lambda_k)^{2m-2}}{\lambda_k
\prod\limits_{i=1,i\neq k\atop i\neq
5-k}^{4}(\lambda_k-\lambda_i)}-
\frac{(1-\lambda_k)^2}{\lambda_k}\right)=
$$
$$
=\frac{2e^{h\beta+h}}{p}\bigg(\frac{(1-\lambda_1)^{4}}{\lambda_1
(\lambda_1-\lambda_2)(\lambda_1-\lambda_3)}-
$$
$$
-\frac{(1-\lambda_1)^2}{\lambda_1}+\frac{(1-\lambda_2)^{4}}{\lambda_2
(\lambda_2-\lambda_1)(\lambda_2-\lambda_4)}-
\frac{(1-\lambda_2)^2}{\lambda_2}\bigg).
$$
Hence, taking into account $\lambda_1\lambda_4=1$,
$\lambda_2\lambda_3=1$ and after some calculations we have
$$
F_1=\frac{2e^{h\beta+h}}{p}\left(\frac{(1-\lambda_1)^2(1-\lambda_2)^2}{
(\lambda_1-\lambda_2)(\lambda_1\lambda_2-1)}-
\frac{(1-\lambda_1)^2(1-\lambda_2)^2}{
(\lambda_1-\lambda_2)(\lambda_1\lambda_2-1)}\right)=0.
$$

\textbf{2).} Now consider convolution
$F_2=D_m(h\beta)*e^{-h\beta}$.
$$
F_2=D_m(h\beta)*e^{-h\beta}=\sum_{\gamma=-\infty}^{\infty}D_m(h\gamma)e^{-h\beta+h\gamma}=
e^{-h\beta}\sum_{\gamma=-\infty}^{\infty}D_m(h\gamma)e^{h\gamma}.
$$
Since $D_m(h\beta)$ is even function, then
$$
F_2=e^{-2h\beta}F_1.
$$
Hence, taking account of 1), when $m=1,2,3$ we obtain, that
$F_2=0$.

\textbf{3).} Consider convolution $F_3=D_m(h\beta)*(h\beta)^n.$

Let $n=0$, then in $m\ge 2$ taking into account (9) we have
$$
D_m(h\beta)*1=\sum_{\gamma=-\infty}^{\infty}D_m(h\gamma)=2\sum_{\gamma=2}^{\infty}
+2D_m(h)+D_m(0)=
$$
$$
=\frac{1}{p}\left(\sum_{k=1}^{m-1}A_k\frac{1+\lambda_k}{\lambda_k(1-\lambda_k)}-
4e^h+2C\right). \eqno(16) $$ For $m=2$ from (16), taking account
of equalities (10), (11), we have
$$
D_2(h\beta)*1=\frac{1}{p}\left(A_1\frac{1+\lambda_1}{\lambda_1(1-\lambda_1)}
-e^h\frac{\lambda_1^2+1}{\lambda_1}+2(e^{2h}+1)\right)=
$$
$$
=\frac{1}{p}\left(-2(e^{2h}+1)+2e^h\frac{\lambda_1^2+1}{\lambda_1}
+2(e^{2h}+1)-2e^h\frac{\lambda_1^2+1}{\lambda_1}\right)=0.
\eqno(17) $$ For $m=3$ from (16), taking account of (10), we get
$$
D_3(h\beta)*1=
$$
$$
=\frac{1}{p}\left(\frac{A_1(1+\lambda_1)}{\lambda_1(1-\lambda_1)}+
\frac{A_2(1+\lambda_2)}{\lambda_2(1-\lambda_2)}-\frac{2e^h(\lambda_1^2+1)}{\lambda_1}
-\frac{2e^h(\lambda_2^2+1)}{\lambda_2}+2(e^h+1)^2\right).
$$
Hence, using (11), after simplifications we obtain
$$
D_3(h\beta)*1=\frac{1}{p}\bigg(-2(e^h+1)^2+2e^h\frac{\lambda_1^2+1}{\lambda_1}
+2e^h\frac{\lambda_2^2+1}{\lambda_2}+
$$
$$ +2(e^h+1)^2-2e^h\frac{\lambda_1^2+1}{\lambda_1}
-2e^h\frac{\lambda_2^2+1}{\lambda_2}\bigg)=0. \eqno(18)
$$

Let now $n=1$, then, using equalities (17), (18) and taking into
account of evenness of the function $D_m(h\beta)$, for $m=2$ we
get$$
D_2(h\beta)*(h\beta)=\sum_{\gamma=-\infty}^{\infty}D_2(h\gamma)(h\beta-h\gamma)=
h\beta\sum_{\gamma=-\infty}^{\infty}D_2(h\gamma)-
\sum_{\gamma=-\infty}^{\infty}D_2(h\gamma)(h\gamma)=0\eqno(19)$$
and also for $m=3$  we obtain $$ D_3(h\beta)*(h\beta)=0.
\eqno(20)$$

Let $n=2$, then taking into account (9), (18), (20) for $m=3$ we
have
$$
D_3(h\beta)*(h\beta)^2=\sum_{\gamma=-\infty}^{\infty}D_3(h\gamma)(h\beta-h\gamma)^2=
\sum_{\gamma=-\infty}^{\infty}D_3(h\gamma)(h\gamma)^2=
$$
$$
=\frac{2h^2}{p}\left(\sum_{k=1}^2\frac{A_k(1+\lambda_k)}{(1-\lambda_k)^3}-2e^h\right).
$$
Hence using (11) and after simplifications we get
$$
D_3(h\beta)*(h\beta)^2=\frac{4h^2}{p}\left(\frac{e^h(\lambda_1-\lambda_2)(\lambda_1\lambda_2-1)}
{(\lambda_1-\lambda_2)(\lambda_1\lambda_2-1)}-e^h\right)=0.
\eqno(21)
$$

For $n=3$ and $m=3$, taking account of (18), (20), (21) and
evenness of the function $D_m(h\beta)$, we obtain
$$
D_3(h\beta)*(h\beta)^3=\sum_{\gamma=-\infty}^{\infty}D_3(h\gamma)(h\beta-h\gamma)^3=
\sum_{\gamma=-\infty}^{\infty}D_3(h\gamma)(h\gamma)^3=0.
$$

\textbf{4).} Consider convolution
$F_4=D_m(h\beta)*\psi_m(h\beta)=\delta(h\beta)$. Taking into
account evenness of functions $D_m(h\beta)$ and $\psi_m(h\beta)$
and definition of convolution of discrete functions, we obtain
$$
D_m(h\beta)*\psi_m(h\beta)|_{\beta=0}=\sum_{\gamma=-\infty}^{\infty}D_m(h\gamma)
\psi_m(h\beta-h\gamma)|_{\beta=0}=\sum_{\gamma=-\infty}^{\infty}D_m(h\gamma)\psi(h\beta).
$$
Hence, taking account of (6) and (9) we have
$$
D_m(h\beta)*\psi_m(h\beta)|_{\beta=0}=\frac{2}{p}\left(-2e^h+\sum_{k=1}^{m-1}A_k\right)
\cdot
\frac{1}{2}\left(\frac{e^h-e^{-h}}{2}-\sum_{k=1}^{m-1}\frac{h^{2k-1}}{(2k-1)!}\right)+
$$
$$
+\frac{2}{p}\sum_{\gamma=2}^{\infty}\sum_{k=1}^{m-1}A_k\lambda_k^{\gamma-1}\frac{1}{2}
\left(\frac{e^{h\gamma}-e^{-h\gamma}}{2}-
\sum_{j=1}^{m-1}\frac{(h\gamma)^{2j-1}}{(2j-1)!}\right)=
$$
$$
=\frac{1}{p}\left(1-e^{2h}+2e^h\sum\limits_{k=1}^{m-1}\frac{h^{2k-1}}{(2k-1)!}\right)+
\frac{1}{p}\sum_{k=1}^{m-1}A_k\left(\frac{e^h-e^{-h}}{2}-\sum_{j=1}^{m-1}
\frac{h^{2j-1}}{(2j-1)!}\right)+
$$
$$
+\frac{1}{p}\sum_{k=1}^{m-1}\frac{A_k}{\lambda_k}\sum_{\gamma=2}^{\infty}\lambda_k^{\gamma}
\left(\frac{e^{h\gamma}-e^{-h\gamma}}{2}-
\sum_{j=1}^{m-1}\frac{(h\gamma)^{2j-1}}{(2j-1)!}\right). \eqno(22)
$$ Using (13) and taking into account convergence of infinite series in equality
(22), we get
$$
D_m(h\beta)*\psi_m(h\beta)|_{\beta=0}=1+\frac{1}{p}\sum_{k=1}^{m-1}\frac{A_k}{\lambda_k}
\bigg(\frac{\lambda_ke^h}{2(1-\lambda_ke^h)}-
$$
$$
-\frac{\lambda_ke^{-h}}{2(1-\lambda_ke^{-h})}-\sum_{j=1}^{m-1}\frac{h^{2j-1}}{(2j-1)!}
\sum_{\gamma=1}^{\infty}\lambda_k^{\gamma}\gamma^{2j-1}\bigg).
\\eqno(23)
$$
Using definition and properties of Euler polynomial, and also (3),
we get
$$
\sum_{\gamma=1}^{\infty}\lambda_k^{\gamma}\gamma^{2j-1}=\frac{1}{1-\lambda_k}\sum_{i=0}
^{2j-1}\left(\frac{\lambda_k}{1-\lambda_k}\right)^i\Delta^i0^{2j-1}=
$$
$$
=\frac{\lambda_k^{2j-1}}{(1-\lambda_k)^{2j}}E_{2j-2}\left(\frac{1}{\lambda_k}\right)=
\frac{\lambda_k}{(1-\lambda_k)^{2j}}E_{2j-2}(\lambda_k). \eqno(24)
$$
Here $E_{2j-2}(\lambda)$ is Euler polynomial of degree $2j-2$.\\
From (23), taking account of (24), we have
$$
D_m(h\beta)*\psi_m(h\beta)|_{\beta=0}=1+\frac{1}{2p}
\sum_{k=1}^{m-1}A_k\bigg(\frac{e^h}{1-\lambda_ke^h}-\frac{1}{e^h-\lambda_k}-
$$
$$
-2\sum_{j=1}^{m-1}\frac{h^{2j-1}E_{2j-2}(\lambda_k)}{(2j-1)!\cdot(1-\lambda_k)^{2j}}\bigg).
$$
Hence, keeping in mind that $\lambda_k$ is a root of polynomial
$P_{2m-2}(\lambda)$ and using (12), we get
$$
D_m(h\beta)*\psi_m(h\beta)|_{\beta=0}=1+\frac{1}{2p}
\sum_{k=1}^{m-1}A_k\frac{P_{2m-2}(\lambda_k)}{(\lambda_ke^h-1)(e^h-\lambda_k)
(1-\lambda_k)^{2m-2}}=1.
$$
Should be noted, that he we obtained
$$
\sum_{\gamma=1}^{\infty}\lambda_k^{\gamma}\psi_m(h\gamma)=
\frac{1}{2}\sum_{\gamma=1}^{\infty}\lambda_k^{\gamma}\left(\frac{e^{h\gamma}-e^{-h\gamma}}{2}-
\sum_{k=1}^{m-1}\frac{(h\gamma)^{2j-1}}{(2j-1)!}\right)=
$$
$$ = \frac{P_{2m-2}(\lambda_k)}{2(\lambda_ke^h-1)(e^h-\lambda_k)
(1-\lambda_k)^{2m-2}}=0 \eqno(25)
$$

Now we prove, that in any $\beta\neq 0$, $\beta$ is a integer and
$m=1,2,3$ following holds
$$
D_m(h\beta)*\psi_m(h\beta)=0
$$
Using the definition of convolution of discrete functions we get
$$
D_m(h\beta)*\psi_m(h\beta)=\sum_{\gamma=1}^{\infty}D_m(h\beta+h\gamma)\psi_m(h\gamma)+
\sum_{\gamma=1}^{\infty}D_m(h\beta-h\gamma)\psi_m(h\gamma).
$$
From here one can see, that it is sufficient to consider case
$\beta>0$.

Let $\beta>0$, then
$$
D_m(h\beta)*\psi_m(h\beta)=\sum_{\gamma=1}^{\beta-2}D_m(h\beta-h\gamma)\psi_m(h\gamma)+
D_m(h)\psi_m(h\beta-h)+
$$
$$
+D_m(0)\psi_m(h\beta)+D_m(-h)\psi_m(h\beta+h)+\sum_{\beta+2}^{\infty}D_m(h\gamma-h\beta)
\psi_m(h\gamma)+
$$
$$
+\sum_{\gamma=1}^{\infty}D_m(h\beta+h\gamma)\psi_m(h\gamma). \eqno
$$
Hence, keeping in mind (9), we obtain
$$
D_m(h\beta)*\psi_m(h\beta)=\frac{1}{p}\sum_{\gamma=1}^{\beta-2}\sum_{k=1}^{m-1}A_k
\lambda_k^{\beta-\gamma-1}\psi_m(h\gamma)+D_m(h)\psi_m(h\beta-h)+
$$
$$
+D_m(0)\psi_m(h\beta)+D_m(-h)\psi_m(h\beta+h)+
\frac{1}{p}\sum_{\beta+2}^{\infty}\sum_{k=1}^{m-1}A_k\lambda_k^{\gamma-\beta-1}
\psi_m(h\beta)+
$$
$$
+\frac{1}{p}\sum_{\gamma=1}^{\infty}\sum_{k=1}^{m-1}A_k\lambda_k^{\beta+\gamma-1}
\psi_m(h\gamma)=D_m(h)\psi_m(h\beta-h)+
$$
$$
+D_m(0)\psi_m(h\beta)+D_m(-h)\psi_m(h\beta+h)+
\frac{1}{p}\sum_{k=1}^{m-1}A_k\lambda_k^{\beta-1}
\sum_{\gamma=1}^{\beta-2}\lambda_k^{-\gamma}\psi_m(h\gamma)-
$$
$$
-\frac{1}{p}\sum_{k=1}^{m-1}A_k\lambda_k^{-\beta-1}
\sum_{\gamma=1}^{\beta+1}\lambda_k^{\gamma}\psi_m(h\gamma)+
\frac{1}{p}\sum_{k=1}^{m-1}A_k\lambda_k^{-\beta-1}
\sum_{\gamma=1}^{\infty}\lambda_k^{\gamma}\psi_m(h\gamma)+
$$
$$ +\frac{1}{p}\sum_{k=1}^{m-1}A_k\lambda_k^{\beta-1}
\sum_{\gamma=1}^{\infty}\lambda_k^{\gamma}\psi_m(h\gamma)
\eqno(26)
$$
We separately calculate
$$
A=\sum_{\gamma=1}^{\beta-2}\lambda_k^{-\gamma}\psi_m(h\gamma)\mbox{
and }
B=\sum_{\gamma=1}^{\beta-2}\lambda_k^{\gamma}\psi_m(h\gamma).
$$
à). Denoting $\lambda_{2k}=\lambda_k^{-1}=\lambda_{2m-1-k}$ and
using (6), we have
$$
A=\sum_{\gamma=1}^{\beta-2}\lambda_k^{-\gamma}\psi_m(h\gamma)=
\sum_{\gamma=1}^{\beta-2}\lambda_k^{-\gamma}\frac{1}{2}
\left(\frac{e^{h\gamma}-e^{-h\gamma}}{2}-\sum_{k=1}^{m-1}
\frac{(h\gamma)^{2j-1}}{(2j-1)!}\right)=
$$
$$
=\frac{1}{4}\left(\sum_{\gamma=1}^{\beta-2}\lambda_{2k}^{\gamma}e^{h\gamma}-
\sum_{\gamma=1}^{\beta-2}\lambda_{2k}^{\gamma}e^{-h\gamma}-2
\sum_{j=1}^{m-1}\frac{h^{2j-1}}{(2j-1)!}
\sum_{\gamma=1}^{\beta-2}\lambda_{2k}^{\gamma}\gamma^{2j-1}\right).
$$
From here, using formulas (2) and (3), we get
$$
A=\frac{1}{4}\Bigg(\frac{\lambda_{2k}e^h}{1-\lambda_{2k}e^h}-\frac{\lambda_{2k}}
{e^h-\lambda_{2k}}-2\sum_{j=1}^{m-1}\frac{h^{2j-1}}{(2j-1)!\cdot(1-\lambda_{2k})}
\sum_{i=0}^{2j-1}\left(\frac{\lambda_{2k}}{1-\lambda_{2k}}\right)^i\Delta^i0^{2j-1}-
$$
$$
-\frac{(\lambda_{2k}e^h)^{\beta-1}}{1-\lambda_{2k}e^h}+\frac{(\lambda_{2k}e^{-h})^{\beta-1}}
{1-\lambda_{2k}e^{-h}}+2\sum_{j=1}^{m-1}\frac{h^{2j-1}\lambda_{2k}^{\beta-1}}
{(2j-1)!\cdot(1-\lambda_{2k})}\sum_{i=0}^{2j-1}\left(\frac{\lambda_{2k}}{1-\lambda_{2k}}\right)^i
\Delta^i(\beta-1)^{2j-1}\Bigg).
$$
Whence, taking into account definition and properties of Euler
polynomial and also formula  (12), we obtain
$$
A=\frac{1}{4}\Bigg(\frac{\lambda_{2k}P_{2m-2}(\lambda_{2k})}{(\lambda_{2k}e^h-1)
(e^h-\lambda_{2k})(1-\lambda_{2k})^{2m-1}}-\frac{(\lambda_{2k}e^h)^{\beta-1}}{1-\lambda_{2k}e^h}
+
$$
$$
+\frac{(\lambda_{2k}e^{-h})^{\beta-1}}
{1-\lambda_{2k}e^{-h}}+2\sum_{j=1}^{m-1}\frac{h^{2j-1}\lambda_{2k}^{\beta-1}}
{(2j-1)!\cdot(1-\lambda_{2k})}\sum_{i=0}^{2j-1}\left(\frac{\lambda_{2k}}{1-\lambda_{2k}}\right)^i
\Delta^i(\beta-1)^{2j-1}\Bigg).
$$
Since $\lambda_k$ is a root of polynomial $P_{2m-2}(\lambda)$,
then we have
$$
A=\frac{1}{4}\Bigg(-\frac{(\lambda_{2k}e^h)^{\beta-1}}{1-\lambda_{2k}e^h}
+\frac{(\lambda_{2k}e^{-h})^{\beta-1}} {1-\lambda_{2k}e^{-h}}+ $$
$$+
2\sum_{j=1}^{m-1}\frac{h^{2j-1}\lambda_{2k}^{\beta-1}}
{(2j-1)!\cdot(1-\lambda_{2k})}\sum_{i=0}^{2j-1}
\left(\frac{\lambda_{2k}}{1-\lambda_{2k}}\right)^i
\Delta^i(\beta-1)^{2j-1}\Bigg).
$$
á). As in case à), by using formulas (6), (2), (3), (12) and
taking into account that $\lambda_k$ is a root of polynomial
$P_{2m-2}(\lambda)$, we have
$$
B=\sum_{\gamma=1}^{\beta-2}\lambda_k^{\gamma}\psi_m(h\gamma)=
\sum_{\gamma=1}^{\beta-2}\lambda_k^{\gamma}\frac{1}{2}
\left(\frac{e^{h\gamma}-e^{-h\gamma}}{2}-\sum_{k=1}^{m-1}
\frac{(h\gamma)^{2j-1}}{(2j-1)!}\right)=
$$
$$
=\frac{1}{4}\Bigg(
-\frac{(\lambda_{k}e^h)^{\beta-1}}{1-\lambda_{k}e^h}+\frac{(\lambda_{k}e^{-h})^{\beta-1}}
{1-\lambda_{k}e^{-h}}+
$$
$$+
2\sum_{j=1}^{m-1}\frac{h^{2j-1}\lambda_{k}^{\beta-1}}
{(2j-1)!\cdot(1-\lambda_{k})}\sum_{i=0}^{2j-1}\left(\frac{\lambda_{k}}{1-\lambda_{k}}\right)^i
\Delta^i(\beta-1)^{2j-1}\Bigg).
$$
Keeping in mind $A$ and $B$ from (26) we get
$$
D_m(h\beta)*\psi_m(h\beta)=
\frac{1}{4p}\sum_{k=1}^{m-1}A_k\lambda_k^{\beta-1}
\Bigg(-\frac{(\lambda_{2k}e^h)^{\beta-1}}{1-\lambda_{2k}e^h}
+\frac{(\lambda_{2k}e^{-h})^{\beta-1}} {1-\lambda_{2k}e^{-h}}+
$$
$$
+2\sum_{j=1}^{m-1}\frac{h^{2j-1}\lambda_{2k}^{\beta-1}}
{(2j-1)!\cdot(1-\lambda_{2k})}\sum_{i=0}^{2j-1}
\left(\frac{\lambda_{2k}}{1-\lambda_{2k}}\right)^i
\Delta^i(\beta-1)^{2j-1}\Bigg)-
$$
$$-
\frac{1}{4p}\sum_{k=1}^{m-1}A_k\lambda_k^{-\beta-1}
\Bigg(-\frac{(\lambda_{k}e^h)^{\beta-1}}{1-\lambda_{k}e^h}+
\frac{(\lambda_{k}e^{-h})^{\beta-1}} {1-\lambda_{k}e^{-h}}+
$$
$$
+2\sum_{j=1}^{m-1}\frac{h^{2j-1}\lambda_{k}^{\beta-1}}
{(2j-1)!\cdot(1-\lambda_{k})}\sum_{i=0}^{2j-1}\left(\frac{
\lambda_{k}}{1-\lambda_{k}}\right)^i
\Delta^i(\beta-1)^{2j-1}\Bigg)-
$$
$$
-\frac{1}{2p}\sum_{k=1}^{m-1}A_k\lambda_k^{-2}\Bigg(
\frac{e^{h\beta-h}-e^{h-h\beta}}{2}-\sum_{j=1}^{m-1}
\frac{(h\beta-h)^{2j-1}}{(2j-1)!}+
\lambda_k\frac{e^{h\beta}-e^{-h\beta}}{2}-
$$
$$
- \lambda_k\sum_{j=1}^{m-1} \frac{(h\beta)^{2j-1}}{(2j-1)!}+
\lambda_k^2\frac{e^{h\beta+h}-e^{-h-h\beta}}{2}-\lambda_k^2\sum_{j=1}^{m-1}
\frac{(h\beta+h)^{2j-1}}{(2j-1)!}\Bigg)+
$$
$$+
\frac{1}{2p}\left(-2e^h+\sum_{k=1}^{m-1}A_k\right)
\Bigg(\frac{e^{h\beta-h}-e^{h-h\beta}}{2}- \sum_{j=1}^{m-1}
\frac{(h\beta-h)^{2j-1}}{(2j-1)!}\Bigg)+
$$
$$
+\frac{1}{2p}\left(2C+\sum_{k=1}^{m-1}\frac{A_k}{\lambda_k}\right)
\left(\frac{e^{h\beta}-e^{-h\beta}}{2}-\sum_{j=1}^{m-1}
\frac{(h\beta)^{2j-1}}{(2j-1)!}\right)+
$$
$$ +\frac{1}{2p}\left(-2e^h+\sum_{k=1}^{m-1}A_k\right)
\left(\frac{e^{h\beta+h}-e^{-h-h\beta}}{2}-\sum_{j=1}^{m-1}
\frac{(h\beta+h)^{2j-1}}{(2j-1)!}\right). \eqno(27)$$

Now in equality (27) we will group coefficients of $e^{h\beta}$,
$e^{-h\beta}$  and $(h\beta)^n$.

First consider coefficients of $e^{h\beta}$. From (27), using (9),
(10), (11), after some calculation we get
$$
\frac{e^{h\beta+h}}{2p}\sum_{k=1}^{m-1}\left(\frac{(1-\lambda_k)^{2m-2}}
{\lambda_k\prod\limits_{i=1,i\neq k\atop i\neq
2m-1-k}^{2m-2}(\lambda_k-\lambda_i)}-\frac{(\lambda_k-1)^2}{
\lambda_k}\right).
$$
Hence for $m=1,2,3$, taking account of (15), we obtain
$$
\frac{e^{h\beta+h}}{2p}\sum_{k=1}^{m-1}\left(\frac{(1-\lambda_k)^{2m-2}}
{\lambda_k\prod\limits_{i=1,i\neq k\atop i\neq
2m-1-k}^{2m-2}(\lambda_k-\lambda_i)}-\frac{(\lambda_k-1)^2}{
\lambda_k}\right)=0.
$$

Now consider coefficients of $e^{-h\beta}$. From (27), taking into
account equalities (9), (10), (11), (15) for $m=1,2,3$, we get
$$
-\frac{e^{-h\beta+h}}{2p}\sum_{k=1}^{m-1}
\left(\frac{(1-\lambda_k)^{2m-2}}
{\lambda_k\prod\limits_{i=1,i\neq k\atop i\neq
2m-1-k}^{2m-2}(\lambda_k-\lambda_i)}-\frac{(\lambda_k-1)^2}{
\lambda_k}\right)=0.
$$
We consider coefficients of $(h\beta)^n$. From equality (27) we
obtain
$$
F_5=\frac{1}{2p}\sum_{j=1}^{m-1}\frac{h^{2j-1}}{(2j-1)!}\Bigg[
\sum_{k=1}^{m-1}A_k\Bigg(\frac{\lambda_k}{\lambda_k-1}\sum_{i=0}^{2j-1}\frac{1}{(\lambda_k-1)^i}
\Delta^i(\beta-1)^{2j-1}-
$$
$$
-\frac{\lambda_k^{-2}}{1-\lambda_k}\sum_{i=0}^{2j-1}
\left(\frac{\lambda_k}{1-\lambda_k}\right)\Delta^i(\beta-1)^{2j-1}+
\frac{(1-\lambda_k^2)(\beta-1)^{2j-1}}{\lambda_k^2}\Bigg)+
$$
$$ +
2e^h(\beta-1)^{2j-1}-2C\beta^{2j-1}+2e^h(\beta-1)^{2j-1}\Bigg].
\eqno(28)
$$
From equality (28) for $m=1$ immediately we obtain, that $F_5=0$.\\
Let $m=2$, then from (28) we have
$$
F_5=A_1\Bigg(\frac{\lambda_1}{\lambda_1-1}(\beta-1+\frac{1}{\lambda_1-1})-
\frac{1}{\lambda_1^2(1-\lambda_1)}(\beta-1+\frac{\lambda_1}{1-\lambda_1})+
$$
$$
+\frac{(1-\lambda_1^2)(\beta-1)}{\lambda_1^2}\Bigg)+\beta(2e^h-2C+2e^h).
$$
Keeping in mind (10) and (11) we obtain
$$
F_5=2\beta\left(\frac{\lambda_1(e^{2h}+1)-e^h(\lambda_1^2+1)}{\lambda_1}-
\frac{\lambda_1(e^{2h}+1)-e^h(\lambda_1^2+1)}{\lambda_1}\right)+
$$
$$
+A_1\left(\frac{\lambda_1(2-\lambda_1)}{(\lambda_1-1)^2}-
\frac{2\lambda_1-1}{\lambda_1^2(1-\lambda_1)^2}-\frac{1-\lambda_1^2}{\lambda_1^2}\right)=0
$$

Let $m=3$, then from  (28), we get
$$
F_5=\frac{1}{2p}\sum_{j=1}^{2}\frac{h^{2j-1}}{(2j-1)!}\Bigg[
\sum_{k=1}^{2}A_k\Bigg(\frac{\lambda_k}{\lambda_k-1}\sum_{i=0}^{2j-1}\frac{1}{(\lambda_k-1)^i}
\Delta^i(\beta-1)^{2j-1}-
$$
$$
-\frac{\lambda_k^{-2}}{1-\lambda_k}\sum_{i=0}^{2j-1}
\left(\frac{\lambda_k}{1-\lambda_k}\right)\Delta^i(\beta-1)^{2j-1}+
\frac{(1-\lambda_k^2)(\beta-1)^{2j-1}}{\lambda_k^2}\Bigg)+
$$
$$ +
2e^h(\beta-1)^{2j-1}-2C\beta^{2j-1}+2e^h(\beta-1)^{2j-1}\Bigg].
\eqno(29)$$ Hence for $j=1$ using (18), we have
$$
\sum_{k=1}^2A_k\Bigg(\frac{\lambda_k}{\lambda_k-1}\left[\beta-1+\frac{1}{\lambda_k-1}\right]-
\frac{\lambda_k^{-2}}{1-\lambda_k}\left[\beta-1+\frac{\lambda_k}{1-\lambda_k}\right]+
\frac{(1-\lambda_k^2)(\beta-1)}{\lambda_k^2}\Bigg)+
$$
$$
+\beta(4e^h-2C)=
\beta\left(\sum_{k=1}^2A_k\frac{\lambda_k+1}{\lambda_k(\lambda_k-1)}+2e^h\sum_{k=1}^2
\frac{\lambda_k^2+1}{\lambda_k}-2(e^h+1)^2\right)+
$$
$$
+\sum_{k=1}^2A_k\left(\frac{\lambda_k(2-\lambda_k)}{(\lambda_k-1)^2}-
\frac{2\lambda_k-1}{\lambda_k^2(1-\lambda_k)^2}-\frac{1-\lambda_k^2}{\lambda_k^2}\right)=0.
$$
From (29) for $j=2$,  grouping by powers of $\beta$ and after some
simplifications we get
$$
\sum_{k=1}^2\Bigg(\frac{\lambda_k}{\lambda_k-1}\sum_{i=0}^3
\frac{1}{(\lambda_k-1)^i}\Delta^i(\beta-1)^3-\frac{1}{\lambda_k^2(1-\lambda_k)}\sum_{i=0}^3
\left(\frac{1}{1-\lambda_k}\right)\Delta^i(\beta-1)^3+
$$
$$
+\frac{1-\lambda_k^2}{\lambda_k^2}(\beta-1)^3\Bigg)+2e^h(\beta-1)^3-2C\beta^3+2e^h(\beta+1)^3=0.
$$
So we proved that in equality (27) for $m=1,2,3$ all coefficients
of $e^{h\beta}$, $e^{-h\beta}$ and $(h\beta)^n$ are zero. Theorem
3 is proved.\\[0.1cm]

\textbf{References}
\begin{enumerate}
\item
Hamming R.W. Numerical methods for scientists and engeneers. -M.:
Nauka, 1968. - 400p.
\item
Sobolev S.L., Vaskevich V.L. Cubature formulas. -Novosibirsk:
Institute of Mathematics SB of RAS, 1996. -484p.
\item
Shadimetov Kh.M. Optimal formulas of approximate integration for
differentiable functions. PhD dissertation, -Novosibirsk, 1983.
\item
Sobolev S.L. Introduction to the Theory of Cubature Formulas .
-M.: Nauka, 1974. - 808p.
\item
Jamalov Z.J. About one problem of Wiener-Hopf arising from
optimization of quadrature formulas. - In book. Boundary value
problems for differential equations . - Tashkent: Fan, 1975. -
pp.129-150.
\item
Jamalov Z.J. About one difference analogue of the operator
$\frac{d^{2m}}{dx^{2m}}$ and its construction. - In book: Direct
and inverse problems for differential equations with partial
derivative and their applications. - Tashkent: Fan. 1978. -
pp.97-108.
\item
Shadimetov Kh.M. Discrete analogue of the operator
$\frac{d^{2m}}{dx^{2m}}$ and its construction. Problems of
computational and applied mathematics. -Tashkent, 1985. V.79. -
pp.22-35.
\item
Shadimetov Kh.M., Hayotov A.R. Construction of discrete analogue
of the differential operator
$\frac{d^{2m}}{dx^{2m}}-\frac{d^{2m-2}}{dx^{2m-2}}$. Uzbek
Mathematical Journal, 2004. \No 2. -pp.85-95.
\item
Kh.M.Shadimetov, A.R.Hayotov. Construction of discrete
ana\-lo\-gue of a dif\-fe\-rential operator. Uzbek Mathematical
Journal, 2003. \No.2. pp.59-69.
\end{enumerate}
\vspace{0.5cm}

\noindent
Kholmat Makhkambaevich Shadimetov\\
Institute of Mathematics and Information Technologies\\
Uzbek Academy of Sciences\\
Tashkent, 100125\\
Uzbekistan\\ [0.2cm]

\noindent
Abdullo Rakhmonovich Hayotov\\
Institute of Mathematics and Information Technologies\\
Uzbek Academy of Sciences\\
Tashkent, 100125\\
Uzbekistan\\
\textit{E-mail:} abdullo\_hayotov@mail.ru, hayotov@mail.ru.

\end{document}